\documentclass[twoside]{amsart}
\usepackage{amsmath,amsthm,amssymb}

\theoremstyle{plain}
\newtheorem{theorem}{Theorem}
\newtheorem{proposition}[theorem]{Proposition}
\newtheorem{lemma}[theorem]{Lemma}
\newtheorem{corollary}[theorem]{Corollary}

\theoremstyle{remark}
\newtheorem{remark}[theorem]{Remark}
\newtheorem{example}[theorem]{Example}
\newtheorem{question}[theorem]{Question}
\newtheorem{acknowledgments}{Acknowledgments}

\newcommand{\llceil}{\left\lceil}
\newcommand{\rrceil}{\right\rceil}
\newcommand{\llfloor}{\left\lfloor}
\newcommand{\rrfloor}{\right\rfloor}
\renewcommand{\k}{{\bf k}}
\newcommand{\K}{{\bf K}}
\newcommand{\Z}{\mathbb{Z}}

\newcommand{\F}{\mathbb{F}}

\newcommand{\supp}{{\rm supp}}
\newcommand{\Hom}{{\rm Hom}}

\newcommand{\NP}{{\rm NP}}

\newcommand{\ord}{ {\rm ord} }
\newcommand{\Jac}{ {\rm Jac} }
\newcommand{\HS}{{\mathcal{HS}}}
\newcommand{\ts}{\tilde s}
\newcommand{\tC}{\tilde C}
\newcommand{\tK}{\tilde K}
\newcommand{\tS}{\tilde S}
\newcommand{\tT}{\tilde T}

\begin{document}

%Topmatter
\title{Families of supersingular curves in characteristic 2}
\author{Jasper Scholten and Hui June Zhu}
\address{
Mathematisch Instituut,
Katholieke Universiteit Nijmegen,
Postbus 9010, 6500 GL Nijmegen.
The Netherlands.
}
\email{scholten@sci.kun.nl}

\address{
Department of mathematics,
University of California,
Berkeley, CA 94720-3840.
The United States.
}
\email{zhu@alum.calberkeley.org}

\date{\today}

\begin{abstract}
This paper determines all normal form of hyperelliptic supersingular
curves of genus $g$ over an algebraically closed field $F$ of
characteristic $2$ for
$1\leq g\leq 8$.  We also show that every hyperelliptic supersingular
curve of genus $9$ over $F$ has an equation $y^2-y=x^{19}+c^8x^9+c^3x$
for some $c\in\overline\F_2$.  Consequently, the paper determines the
dimensions of the open locus of hyperelliptic supersingular curves of
genus $g\leq 9$ over $\overline\F_2$.
\end{abstract}

\maketitle

\section{Introduction}

In this paper, a curve is a projective, smooth and geometrically
integral algebraic variety of dimension $1$.  A curve is {\em
supersingular} if its Jacobian is a supersingular abelian variety, that
is, its Newton polygon is a straight line segment of slope $1/2$.  See
Introduction and Appendix of \cite{Li-Oort:98} for literatures on
supersingular abelian varieties and related open questions on
stratification in the moduli space of abelian varieties.  It has been
proved that there are supersingular curves of every genus over
$\overline\F_2$ (see \cite{Geer&Vlugt}), but it is not settled if
there are hyperelliptic supersingular curves over $\overline\F_2$ of
every genus, which was an initial goal of our study (see \cite{SZ:1}).
Several isolated discoveries of supersingular curves have yielded
unprecedented applications in sphere packing \cite{Elkies:94} and
\cite{Elkies:97}. There are also rich literatures in coding theorem
upon applications of supersingular curves (see \cite{Geer}).  See the
introduction in \cite{SZ:1} for more recent progress regarding the
hyperelliptic supersingular curves in characteristic $2$.  In
\cite{deJong-Oort} a question is raised that for which $p$ and $g$
there are positive dimensional families of hyperelliptic supersingular
curves over $\overline\F_p$ of genus $g$. We give an answer to this
question for $p=2$ and $g\leq 9$ in Theorem \ref{T:1}.

Every genus-$g$ hyperelliptic curve $X$ over $F$ of $2$-rank
zero has an affine equation
\begin{eqnarray}\label{E:1}
X:\ y^2-y&=&c_{2g+1}x^{2g+1}+c_{2g-1}x^{2g-1}+\ldots+c_1x,
\end{eqnarray}
where $c_1,\ldots,c_{2g-1}\in F$ and $c_{2g+1}=1$.  (See
Proposition \ref{P:1}).

For any $g$, an equation (\ref{E:1})
with some vanishing conditions among $c_\ell$'s
is called a {\em normal form of hyperelliptic
supersingular curve of genus $g$} if every such equation defines a
hyperelliptic supersingular curve and conversely every hyperelliptic
supersingular curve of genus $g$ has such an equation for some
$c_\ell$'s in $F$, and there are only finitely many such
equations for each curve.

Let $\HS_g/F$ denote the open locus of hyperelliptic
supersingular curves of genus $g$ over $F$.

\begin{theorem}\label{T:1}
Let $F$ be an algebraically closed field of characteristic $2$.
The normal forms of
hyperelliptic supersingular curves over $F$ of genus $g\leq
8$ are in the table below with all $c_\ell\in F$.
The dimensions of $\HS_g/F$ are in the right columns.%(Insert Table 1 here.)
\begin{table}[h]
	\begin{center}
\begin{tabular}{|c|l|c|}\hline
genus g&normal form $X_g$&$\dim(\HS_g/F)$\\\hline
$1$&$y^2-y=x^3$&$0$\\\hline
$2$&$y^2-y=x^5+c_3x^3$&$1$\\\hline
$3$&none&$-\infty$\\\hline
$4$&$y^2-y=x^9+c_5x^5+c_3x^3$&$2$\\\hline
$5$&$y^2-y=x^{11}+c_3x^3+c_1x$&$2$\\\hline
$6$&$y^2-y=x^{13}+c_3x^3+c_1x$&$2$\\\hline
$7$&none&$-\infty$\\\hline
$8$&$y^2-y=x^{17}+c_9x^9+c_5x^5+c_3x^3$&$3$\\\hline
\end{tabular}
\end{center}
\caption{normal forms for $g\leq 8$}
\end{table}

Moreover, every supersingular curve over $F$
of genus $9$ has an equation
$$y^2-y=x^{19}+c^8x^9+c^3x$$
for some $c\in\overline{\F}_2$. We have
$\dim(\HS_9/F)=0$.
\end{theorem}

\begin{remark}\label{R:1}
{}From Theorem~\ref{T:1} one can easily derive the following:
If $g=1$ or $2$, then all 2-rank zero curves are
supersingular.  If $g=3$ or $7$, then none of these curves is
supersingular. If $g=4$ then $X$ is supersingular if and only if
$c_7=0$.  If $g=5$ then $X$ is supersingular if and only if $c_7=0$
and $c_9c_{11}=c_5^4$. If $g=6$ then $X$ is supersingular if and only
if $c_7=0$, $c_{11}=0$, and $c_5c_{13} = c_9^2$.  If $g=8$ then $X$ is
supersingular if and only if $c_{15}=c_{13}=c_{11}=c_{7}=0$.
\end{remark}

\begin{remark}
It is clear from Theorem \ref{T:1} that
$\HS_g/F$ for $1\leq g\leq 8$ are irreducible and unirational.
One wonders if it is so for all $g$.
\end{remark}

The proof of Theorem \ref{T:1} splits into two parts:
In section \ref{SS:supersingular} we show that all curves
in Table 1 are supersingular. In section \ref{S:necessarity}
we show that every supersingular curve of genus $\leq 9$
has an equation in Table 1.

One main idea in the proof of Theorem \ref{T:1} is to use a lemma of
\cite{SZ:2} to reduce the initial problem to a combinatorial puzzle.
Then we formulate an algorithm, utilizing Proposition \ref{P:1} and
Theorem \ref{MT:3} from \cite{SZ:2}, to prove our assertion.

We anticipate some results of similar sort for Artin-Schreier curves
in characteristic $p$ if $p$ is small compared to the genus.
For some results for $p$ large, see \cite{Zhu:1}.

\begin{acknowledgments}
We are grateful to Noam Elkies for his suggestion
in section \ref{SS:supersingular}.
We especially thank Bjorn Poonen for his valuable suggestions
in Remark \ref{R:1}.
The computation in this paper used Magma and Maple
packages and C-programs.
\end{acknowledgments}

\section{$2$-adic box and a tiling puzzle}\label{SS:puzzle}

Let $s(\cdot)$ denote the number of $1$'s in
the binary expansion of a positive integer.
This section digresses to some new combinatorial notion.
Let $S$ be a finite set of positive integers.
Let $r$ be a positive integer.
For a finite sequence of integer pairs
$\{[b_i,\ell_i]\}_{i=1}^{\nu}$ such that
\begin{enumerate}
\item $\ell_i\in S$, $0\leq b_i\leq b_{i+1}$;
\item $\ell_i>\ell_{i+1}$ if $b_i=b_{i+1}$;
\item $\sum_{i=1}^\nu\ell_i2^{b_i}=r$
\end{enumerate}
we call it {\em a tiling sequence of $r$} or an {\em $r$-tiling
sequence}.  If no such sequence exists we set $\ts(r,S):=\infty$;
otherwise, we denote by $\ts(r,S)$ the minimal length among all tiling
sequences of $r$. Let $\tK(r,S)$ denote the set of all such minimal
length tiling sequences of $r$. For any $r$-tiling sequence
$M=\{[b_i,\ell_i]\}_{i=1}^{\nu}$, let $\tS(M):=\{\ell_1,\ldots,\ell_\nu\}$.

\begin{example}\label{Ex:1}
Put $r=2^{17}-2^2$ and $S=\{1,3,9,11\}$.
Observe that the binary expansion of $r$ has $15$ consecutive
$1$'s followed by two $0$'s. To make the shortest $r$-tiling
sequence ideally one has to choose $\ell\in S$
with highest $s(\ell)$ possible and with every $1$'s contributes to
$r$ when added up with each other. In this example,
$(9)_2=1001$ and $(11)_2=1011$ make $5$ consecutive $1$'s if added as
$$
\begin{array}{cccccc}
   &  &1&0&0&1\\
+) & 1&0&1&1& \\
               \hline
   & 1&1&1&1&1
\end{array}
$$
Thus to make $15$ consecutive $1$'s we just have to combine
$3$ such blocks. The minimal length $r$-tiling sequence
is demonstrated below.
\begin{table}[h]
\begin{center}
$$
\begin{array}{cccccccccccccccccc|l} %6|5|5|2
& & & & & & & & & & & &1&0&0&1&0&0&[2,9]\\
& & & & & & & & & & &1&0&1&1&0&0&0&[3,11]]\\
& & & & & & &1&0&0&1&0&0&0&0&0&0&0&[7,9]\\
& & & & & &1&0&1&1&0&0&0&0&0&0&0&0&[8,11]\\
& &1&0&0&1&0&0&0&0&0&0&0&0&0&0&0&0&[12,9]\\+)
&1&0&1&1&0&0&0&0&0&0&0&0&0&0&0&0&0&[13,11]\\
\hline\hline
&1&1&1&1&1&1&1&1&1&1&1&1&1&1&1&0&0
\end{array}
$$
\end{center}
\caption{The minimal length $(2^{17}-2^2)-$ tiling}
\end{table}
The only way to make a shorter $r$-tiling sequence is to
use more of the form $[*,11]$'s in the sequence, which is
impossible since more $11$ will not make consecutive $1$'s
in its sum unless some $1$'s does not contribute to the sum.
So $\ts(r,S)=6$ and
$\tK(r,S)=\{\{[2,9],[3,11],[7,9],[8,11],[12,9],[13,11]\}\}$.
It is computationally effective to represent the set $\tK(r,S)$
binarily as in the above table. But since it takes too much space,
in this paper we shall adopt the following way of representation
\begin{eqnarray*}
2^{17}-2^2
&=& 11\cdot 2^{13} + 9\cdot 2^{12} + 11\cdot 2^8 + 9\cdot 2^7
 +11\cdot 2^3  + 9\cdot 2^2.
\end{eqnarray*}
Denote this minimal length $r$-tiling sequence by $M$ then
$\tS(M)=\{9,11\}$.
\end{example}

Set $d:=2g+1$ for the rest of the paper.  Recall the $2$-adic box
$\fbox{\k}$ and notations defined in \cite[Section 2]{SZ:1}.
\begin{eqnarray*}
\K_r&:=&\{\k={}^t\!(k_1,k_2,\ldots,k_d)\in\Z^d\mid
k_1\geq k_2\geq\cdots\geq k_d\geq 0, \sum_{\ell=1}^{d} k_\ell = r
\}\\
s(\k)&=& s(k_1-k_2)+s(k_2-k_3)+\ldots + s(k_{d-1}-k_{d}) + s(k_{d}).
\end{eqnarray*}

\begin{lemma}\label{L:bijection}
There is a bijection between the sets
$\tK(r,S)$ and $$
\{\k\in\K_r\ |\ s(\k)=\ts(r,S)\ {\rm and}\ k_\ell=k_{\ell+1}
\ {\rm for\ all}\ \ell\not\in S\}.$$
\end{lemma}
\begin{proof}
We shall define the map first.
An $r$-tiling sequence  $\{[b_i,\ell_i]\}_{i=1}^{\ts(r,S)}\in\tK(r,S)$
is sent to the element $\k\in\K_r$ whose {\em 2-adic box} $\fbox{\k}$ has
$ k_{\ell,v}=\#\{i\ |\ v=b_i{\rm\ and\ }\ell\leq\ell_i\}$.

Conversely, given $\k\in\K_r$ with $k_\ell=k_{\ell+1}$ for
$\ell\not\in\supp X$, one defines $\{[b_i,\ell_i]\}_{i=1}^{s(\k)}$ as
follows: For $v$ such that $ k_{1,v}\neq 0$ and for $i$ such
that $\sum_{j=0}^{v-1} k_{1,j}<i\leq\sum_{j=0}^v k_{1,j}$ let
$b_i:=v$. Let $m_1>\ldots>m_{ k_{1,v}}>0$ be the sequence
of positive integers such that $ k_{m_1,v}> k_{m_2,v}
\ldots> k_{m_{ k_{1,v}},v}$ and such that
$ k_{m_j}> k_{m_j+1}$. Note that $m_j\in S$.
Define $\ell_i:=m_{i-\sum_{j=0}^{v-1} k_{j,v}}$.
The detailed verification of the bijectivity of the map is
tedious and we omit this.
Instead we shall visualize the bijection below for Example
\ref{Ex:1}.
\end{proof}

In the Table 3 below, the last line is the binary expansion of
$r=2^{17}-2^2$. The block above the double horizontal lines is the
$2$-adic box $\fbox{\k}$. The $6$ lines immediately beneath the
double horizontal lines are the column sums of $\fbox{\k}$ added up
in the order from right to left.
They represent the minimal length $r$-tiling sequence
of $\tK(r,S)$ as seen in Table 2. Their sums are of course equal to
$r$ as shown below the single line.
\begin{table}[h]
\begin{center}
$$
\begin{array}{cccccccccccccccccc} % 6|5|5|2
&0&0&0&1&1&0&0&0&1&1&0&0&0&1&1&0&0\\ %18
&0&0&0&1&1&0&0&0&1&1&0&0&0&1&1&0&0\\
&0&0&0&1&1&0&0&0&1&1&0&0&0&1&1&0&0\\
&0&0&0&1&1&0&0&0&1&1&0&0&0&1&1&0&0\\
&0&0&0&1&1&0&0&0&1&1&0&0&0&1&1&0&0\\
&0&0&0&1&1&0&0&0&1&1&0&0&0&1&1&0&0\\
&0&0&0&1&1&0&0&0&1&1&0&0&0&1&1&0&0\\
&0&0&0&1&1&0&0&0&1&1&0&0&0&1&1&0&0\\
&0&0&0&1&1&0&0&0&1&1&0&0&0&1&1&0&0\\
&0&0&0&1&0&0&0&0&1&0&0&0&0&1&0&0&0\\+)
&0&0&0&1&0&0&0&0&1&0&0&0&0&1&0&0&0\\
\hline\hline
& & & & & & & & & & & &1&0&0&1&0&0\\
& & & & & & & & & & &1&0&1&1&0&0&0\\
& & & & & & &1&0&0&1&0&0&0&0&0&0&0\\
& & & & & &1&0&1&1&0&0&0&0&0&0&0&0\\
& &1&0&0&1&0&0&0&0&0&0&0&0&0&0&0&0\\+)
&1&0&1&1&0&0&0&0&0&0&0&0&0&0&0&0&0\\
\hline
&1&1&1&1&1&1&1&1&1&1&1&1&1&1&1&0&0
\end{array}
$$
\end{center}
\caption{The $2$-adic box and the minimal length $r$-tiling sequence}
\end{table}

\section{Supersingularity criterion}

Let $X$ be as defined with affine equation in (\ref{E:1}). Let
$\NP_1(X)$ be the first slope of the Newton polygon of $X$ (see
\cite{SZ:1}). Recall that $d:=2g+1$.

Define the set of {\em supports} of $X$ by $\supp X:=\{\ell\in\Z \mid
c_\ell\neq 0\}$.  For $r\geq 1$ let
\begin{eqnarray*}
\tC(r,\supp X)&:=&\sum_{M\in\tK(r,\supp X)}
     \prod_{[b,\ell]\in M}c_{\ell}^{2^b}\in F.
\end{eqnarray*}

\begin{lemma}\label{L:key}
Let notations be as above.
Let $\lambda$ be a rational number with $0\leq \lambda\leq 1/2$.
\begin{enumerate}
\item[a)]
Suppose that for all $m\geq 1$, $n\geq 1$, and  $1\leq j\leq g$
one has
\begin{eqnarray*}
\ts(m2^{n+g-1}-j,\supp X)&\geq&\lceil n\lambda\rceil.
\end{eqnarray*}
Then
\begin{eqnarray*}
\NP_1(X)&\geq&\lambda.
\end{eqnarray*}
\item[b)]
Suppose $\NP_1(X)\geq\lambda$.
Suppose there are positive integers $1\leq j\leq g$, $n_0\geq 1$ and
$1\leq m_0\leq g/2$ such that
\begin{enumerate}
\item[1)] for all $m\geq 1,1\leq n<n_0$ and for all $m>g/2, n=n_0$
we have \begin{eqnarray*}
\ts(m2^{n+g-1}-j, \supp X)&\geq& \llceil n\lambda\rrceil;
\end{eqnarray*}
\item[2)] $\ts(m_02^{n_0+g-1}-j, \supp X) = \lceil n_0\lambda\rceil-1.$
\end{enumerate}
Then
\begin{eqnarray*}
\tC(m_02^{n_0+g-1}-j, \supp X)&=&0.
\end{eqnarray*}
\end{enumerate}
\end{lemma}

\begin{proof}
Let $W(F)$ be the ring of Witt vectors of $F$.
Let $a_\ell\in W(F)$ be a lift of $c_\ell\in F$.
Let
$B:=\min\{s(\k)\ |\ \k\in\K_r,\ {\rm and}\ k_\ell=k_{\ell+1}
\ {\rm for\ all}\ \ell\not\in \supp X\}.$
Set $0^0:=1$. We recall $C_r(N)$ from
\cite{SZ:1}; we know from there that
\begin{eqnarray*}
C_r(N)&\equiv&\sum_{\k\in \K_r}2^{s(\k)}
\prod_{\ell=1}^{d-1}a_\ell^{k_{\ell}-k_{\ell+1}}\bmod 2^{B+1}.
\end{eqnarray*}

Since $\ord_2(C_r(N))\geq B$, we may write
\begin{eqnarray*}%\label{E:C_r}
\frac{C_r(N)}{2^B}&\equiv& \sum_{\k\in\K_r} 2^{s(\k)-B}\prod_{\ell=1}^{d-1}
c_\ell^{k_\ell-k_{\ell+1}} \bmod 2.
\end{eqnarray*}
By the bijection discussed in section \ref{SS:puzzle}
we get
\begin{eqnarray}\label{E:C}
(\frac{C_r(N)}{2^{\ts(r,\supp X)}}\bmod 2)=\tC(r,\supp X).
\end{eqnarray}
By (\ref{E:C}), we note that a) and b) follow from
the Key Lemma of \cite{SZ:2}.
(see also \cite[ Key-Lemma~2.1,i)]{SZ:1}).
\end{proof}

Below we have two corollaries.

\begin{corollary}\label{C:1}
Let $X$ be as in (\ref{E:1}).
Let $\sigma:=\max_{\ell\in\supp X}s(\ell)$. Then we have $\NP_1(X)\geq
1/\sigma\geq 1/d$.
\end{corollary}
\begin{proof}
For $r=m2^{n+g-1}-j$ consider any $r$-tiling sequence
given by $m2^{n+g-1}-j=\sum_{i=1}^\nu\ell_i2^{b_i}$. Clearly we have
$$\nu\sigma
\geq \sum_{i=1}^{\nu}s(\ell_i)
=    \sum_{i=1}^{\nu}s(\ell_i2^{b_1})\geq s(m2^{n+g-1}-j)\geq n,
$$
so $\nu\geq\lceil\frac{n}{\sigma}\rceil.$ In particular
$$\ts(m2^{n+g-1}-j,\supp X)\geq \lceil\frac{n}{\sigma}\rceil.$$
By Lemma~\ref{L:key} a) we have $\lambda\geq 1/\sigma$.
But $c_d=1\neq 0$ so  $\sigma\leq d$ and
we have $\NP_1(X)\geq 1/\sigma\geq 1/d$.
\end{proof}

\begin{remark}\label{R:X_8}
Any curve $X$ given by an equation of the form
$y^2-y=\sum_{i=0}^{s}c_ix^{2^i+1}$ over $F$ is
supersingular. Indeed, apply the above corollary to this curve
we get
$\sigma=2$ and so $\NP_1(X)\geq 1/2$. (This result was proved in
\cite{Geer&Vlugt} in the case $F\subseteq\overline\F_2$.)
\end{remark}

\begin{corollary}\label{C:2}
Suppose $X$ given by equation (\ref{E:1}) is supersingular.  Let
$m_0$, $n_0$ and $j$ be positive integers satisfying hypothesis of
Lemma~\ref{L:key} b) for $\lambda=\frac{1}{2}$.  Suppose the set
$\tK(m_02^{n_0+g-1}-j, \supp X)$ consists of one element $M$.  Then we
have $\tS(M)\neq\{2g+1\}$ and if $\tS(M)=\{\ell,2g+1\}$ with
$1\leq\ell\leq 2g$, then $c_{\ell}=0$.
\end{corollary}
\begin{proof}
Suppose $\tS(M)=\{2g+1\}$, then
by Lemma \ref{L:key} b) we have $\tC(m_02^{n_0+g-1}-j,\supp X)=0$.
But it is easily seen that $$\tC(m_02^{n_0+g-1}-j,\supp X)
=\prod_{[b,2g+1]\in M}c_{2g+1}^{2^b}=1$$
because $c_{2g+1}=1$.
Contradiction. This proves the first assertion.

Now suppose $\tS(M)=\{\ell,2g+1\}$.
By Lemma \ref{L:key} again, we have
\begin{eqnarray*}
\tC(m_02^{n_0+g-1}-j,\supp X)
&=&\prod_{[b,\ell]\in M}c_\ell^{2^{b}}\prod_{[b',2g+1]\in M}c_{2g+1}^{2^{b'}}\\
&=&\prod_{[b,\ell]\in M}c_\ell^{2^{b}}
=c_\ell^{\sum_{[b,\ell]\in M}2^b}=0.
\end{eqnarray*}
Thus $c_\ell=0$.
\end{proof}

\section{Supersingularity}\label{SS:supersingular}

We shall show that all curves in Table 1 are supersingular.

The supersingularity of $X_8$ follows from Remark \ref{R:X_8}.  The
supersingularity of the curve $C: y^2-y = x^{33}+c_3x^9+c_1x^3$ over
$F$ also follows from Remark \ref{R:X_8}. By observing that $C$
covers $X_5$,  the supersingularity of $X_5$ follows.

Now we claim that $X_6$ is supersingular
by a method suggested to us by Noam Elkies. We present it below:
Let $E$ be the elliptic curve defined by
$y^2-y=x^3+t^{13}+c_3t^3+c_1t$ over $F(t)$, which can be viewed
as a quadratic twist of $E_0: y^2-y = x^3$ over the function field
$F(X_6)$ of $X_6$.  The curve $X_6$ is supersingular if and only if
the rank of $\Hom_{F}(\Jac(X_6),E_0)$, or equivalently the rank of
$E(F(t))$ is $24$, the maximal possible
(see~\cite{Sha-Tate}). In the same vein as in~\cite[pages
3--4]{Elkies:97}, one can show that the canonical height of a nonzero
point in $E(F(t))$ is $\geq 8$. We shall show in the next
paragraph that there
are $196560$ nonzero points with the minimal canonical height
$8$. Since there are no lattices of rank $\leq 23$ with that many
vectors of minimal length~\cite[page 23]{CS:93}, this implies that
$E(F(t))$ has rank $24$ and the supersingularity of $X_6$
follows.

If $P=(x,y)$ is a point in $E(F(t))$ with $x$-coordinate
of the form
$$x=a^{-16}t^6+a^{-28}t^4+a^{-40}t^2+a^{32}t+x_0+a^{-64}(t-t_0)^{-2},$$
for $a\in F^*$ and $x_0,t_0\in F$, then $P$ has canonical
height $8$ (see~\cite[Proposition~2]{Elkies:94}).
Following closely the computations presented in~\cite[page 6--7]{Elkies:97},
one finds that $P$ lies in $E(F(t))$
if $a,x_0$ and $t_0$ satisfy the following equations:
\begin{eqnarray}
\label{E:inf1}
a^{4096}+c_3^{32}a^{1024}+c_1^{16}a^{256}+c_3^{16}a^{64}+
c_1^4a^{16}+c_3^2a^4+a &=& 0;\\
\label{E:inf2}
a^{36}t_0^6+a^{24}t_0^4+a^{12}t_0^2+a^{84}t_0+
a^{1092}+c_3^8a^{324}+&&\\ \nonumber
a^{312}+a^{156}+c_1^4a^{132}+c_3^2a^{120}+a^{117}+c_3a^{84}+1 &=&0;\\
\label{E:inf3}
{a}^{576}+{a}^{264}+{c_3}^{4}{a}^{192}+
(x_0^{4}+x_0){a}^{160}+{a}^{108}+t_0^{8} &=&0.
\end{eqnarray}
(One may also check alternatively that the solutions to these $3$
equations indeed yield $x$-coordinates of points in $E(F(t))$.)
One can show that these $3$ equations are separable in $a,t_0$ and
$x_0$ respectively by observing that their derivatives are nonzero
constant. The number of nonzero solutions $a$ of $(\ref{E:inf1})$ is
$4095$. For each such $a$ there are $6$ solutions of $t_0$ satisfying
$(\ref{E:inf2})$.  For each pair of $(a,t_0)$ there are $4$ solutions
of $x_0$ satisfying $(\ref{E:inf3})$. So there are totally $4095\cdot
6\cdot 4$ possible $x$-coordinates, each of which yields two points.
This results in $2\cdot 4\cdot 6\cdot 4095=196560$ points in
$E(F(t))$.

\section{Necessary}\label{S:necessarity}

Lemma~\ref{L:key} b) shall find its extensive application here.  A few
remarks are in order. If there are indeed positive integers $j, n_0,
m_0$ satisfying the hypothesis of Lemma~\ref{L:key} b) for some
$\lambda$, then it is a finite problem to find them and verify.  To
see this, define $\tT(\nu,j,\supp X)$ as the set of all sequences
$\{[b_i,\ell_i]\}_{i=1}^v$ satisfying b) and that
$\sum_{i=1}^\nu\ell_i2^{b_i}\equiv -j\bmod 2^{b_\nu}$.  Clearly,
$\tT(\nu,j,\supp X)$ is a finite set, and if
$\nu=\ts(m_02^{n_0+g-1}-j,\supp X)$ for some $m_0$ and $n_0$, then
$\tT(\nu,j,\supp X)$ contains $\tK(m2^{n_0+g-1}-j,\supp X)$ for every
$m$ with $$\ts(m2^{n_0+g-1}-j,\supp X)\leq \ts(m_02^{n_0+g-1}-j,\supp
X).$$  So one can verify if b) is satisfied by computing $\tT(\nu,j,\supp
X)$ for $1\leq j\leq g$ and for $\nu\geq 1$, and check for each such
$\nu$ and $j$ whether there exist $m_0$ and $n_0$ with the required
property.

We shall verify the hypotheses of b) by the
aforementioned algorithm, together with a few isolated tricks to
accelerate the program. For instance, we use an observation that if
two sequences $\{[b_i,\ell_i]\}_{i=1}^\nu$ and
$\{[b'_i,\ell'_i]\}_{i=1}^{\nu'}$ satisfy
$\sum_{i=1}^\nu\ell_i2^{b_i}=\sum_{i=1}^{\nu'}\ell'_i2^{b'_i}$ and
$\nu'<\nu$, then one can ignore all sequences including
$\{[b_i,\ell_i]\}_{i=1}^\nu$ as a subsequence.

For convenience of the reader we first quote
Proposition~4.1 and Theorem 1.1 of \cite{SZ:1} below.

\begin{proposition}\label{P:1}
Let $X$ be a hyperelliptic curve over $F$ of genus $g$ and
$2$-rank zero. Let $m< 2g$ be an odd integer, such that the binomial
coefficient $\binom{2g+1}{2^km}$ is odd for some $k\geq 0$.  Then
there exists an equation for $X$ of the form
\begin{eqnarray}\label{E:oddterms}
y^2-y=\sum_{i=0}^{g}c_{2i+1}x^{2i+1}
\end{eqnarray}
with $c_{2g+1}=1$ and $c_m=0$. Moreover, there are only finitely
many equations for $X$ of this form.
\end{proposition}

\begin{theorem}\label{MT:3}
Let $g\geq 3$. Write $h=\llfloor\log_2(g+1)+1\rrfloor$.
Let $X$ be a hyperelliptic curve
over $F$ given by an equation $y^2 - y = \sum_{i=1}^{2g+1}c_ix^i$
where $c_{2g+1}=1$.  Then
\begin{enumerate}
\item[(I)] $\NP_1(X)\geq\frac{1}{h}$.
\item[(II)]
If $g<2^h-2$ and $c_{2^{h}-1}\neq 0$ then
$\NP_1(X)=\frac{1}{h}$.
\item[(III)]
If $g=2^h-2$ and
$c_{2^{h}-1}\neq 0$ or $c_{3\cdot2^{h-1}-1}\neq 0$ then
$\NP_1(X) = \frac{1}{h}$.
\end{enumerate}
\end{theorem}

Now we will proceed with the proof of Theorem~\ref{T:1} for genus $g\geq 5$.
For the rest of the proof we shall suppress the $\supp X$ from the notations
$\tC(r,\supp X)$ and $\tK(r,\supp X)$. For every genus we start
with a curve with an equation in a box. Then we use
supersingularity criterion to derive vanishness of coefficients
of this initial equation, which accumulate to our desired equation
declared in Table 1.

\subsection{$g=5$}

Suppose $X$ is a hyperelliptic supersingular curve over $F$ of genus $5$.
By Theorem~\ref{MT:3} II) and Proposition~\ref{P:1}
it has an equation
$$\fbox{$y^2-y=x^{11}+c_9x^9+c_3x^3+c_1x$}$$
for some $c_1,c_3,c_9\in F$.

Let $ \supp X=\{1,3,9,11\}$.
Let $g=5$. We checked that for $n_0=13$, $m_0=1$ and $j=4$ the
conditions of Lemma~\ref{L:key} b) are satisfied, and
$\tK(m_02^{n_0+g-1}-j)$ consists of a unique element $M$ as follows
$$
2^{17}-4
= 11\cdot 2^{13} + 9\cdot 2^{12} + 11\cdot 2^8 + 9\cdot 2^7
 +11\cdot 2^3  + 9\cdot 2^2.
$$
Note that $\tS(M)=\{9,11\}$.
(See Example \ref{Ex:1}.)
Then Corollary \ref{C:2} implies that \fbox{$c_9=0$}.

\subsection{$g=6$}

Suppose $X$ is a hyperelliptic
supersingular curve over $F$ of genus $6$.
By Theorem~\ref{MT:3} III) and Proposition~\ref{P:1}
it has an equation
$$\fbox{$y^2-y=x^{13}+c_9x^9+c_3x^3+c_1x$}$$
for some $c_1,c_3,c_9\in F$.

Let $ \supp X=\{1,3,9,13\}$.
Let $g=6$.
We checked that for $n_0=17$, $m_0=1$, $j=4$
the conditions of Lemma~\ref{L:key} b) are satisfied, and
$\tK(m_02^{n_0+g-1}-4)$ consists of
a unique element $M$ as follows
$$2^{22}-4 = 9\cdot 2^{18}+13\cdot 2^{17} +9\cdot 2^{13}+13\cdot
2^{12}+ 9\cdot 2^8+13\cdot 2^7 +9\cdot 2^3+
13\cdot 2^2.$$ Note that $\tS(M)=\{9,13\}$.
By Corollary \ref{C:2} we have \fbox{$c_9=0$}.

\subsection{$g=8$}

Suppose $X$ is a  hyperelliptic supersingular curve over
$F$ of genus $8$. By Theorem~\ref{MT:3} II) and
Proposition~\ref{P:1}
 it has an equation
$$\fbox{$y^2-y=x^{17}+c_{13}x^{13}+c_{11}x^{11}+c_9x^9+c_7x^7+
c_5x^5+c_3x^3$}$$
for some $c_\ell\in F$.

Let $\supp X=\{3,5,7,9,11,13,17\}$. We found for
$k=0, 1, 2, 3, 4$ and 5 that
$\tK(2^{21-3k}-8)$ consists of
sequences $\{[b_i,\ell_i]\}_{i=1}^{6-k}$ with either
$[b_{6-k},\ell_{6-k}]=[18-3k,7]$ or
$[b_{6-k},\ell_{6-k}]=[17-3k,13]$ and
$[b_{5-k},\ell_{5-k}]=[15-3k,11]$.
Hence we get
\begin{eqnarray*}
\tC(2^{21}-8 )&=&c_7^{2^{18}}\tC(2^{18}-8 )+
(c_{11}c_{13}^4)^{2^{15}}\tC(2^{15}-8 )\\
\tC(2^{18}-8 )&=&c_7^{2^{15}}\tC(2^{15}-8 )+
(c_{11}c_{13}^4)^{2^{15}}\tC(2^{12}-8 )\\
\tC(2^{15}-8 )&=&c_7^{2^{12}}\tC(2^{12}-8 )+
(c_{11}c_{13}^4)^{2^{15}}\tC(2^{9}-8 )\\
\tC(2^{12}-8 )&=&c_7^{2^{9}}\tC(2^{9}-8 )+
(c_{11}c_{13}^4)^{2^{15}}\tC(2^{6}-8 )\\
\tC(2^{9}-8 )&=&c_7^{2^{6}}\tC(2^{6}-8 )+
(c_{11}c_{13}^4)^{2^{15}}\tC(2^{3}-8 )\\
\tC(2^{6}-8 )&=&c_7^{2^{3}}\tC(2^{3}-8 ).
\end{eqnarray*}
We checked that the conditions of Lemma~\ref{L:key} b)
are satisfied for
$\lambda=\frac{1}{2}$, $n_0=11$, $m_0=1$ and $j=8$. Hence
$\tC(2^{18}-8 )=0$.
We also checked that they are satisfied for
$\lambda=\frac{4}{9}$, $n_0=14$, $m_0=1$ and $j=8$. Hence
$\tC(2^{21}-8 )=0$.
All these equations imply \fbox{$c_7=0$} and $c_{11}c_{13}=0$.

Suppose $c_{11}=0$. Let $\supp X=\{3,5,9,13,17\}$. We checked that
for $\lambda=\frac{1}{2}$, $j=4$, $n_0=19$ and $m_0=1$ the conditions
of Lemma~\ref{L:key} b) are satisfied, and
$\tK(m_02^{n_0+g-1}-4 )$ consists of a unique element $M$ as
follows
\begin{eqnarray*}
2^{26}-4&=&13\cdot2^{22}+17\cdot2^{19}+13\cdot2^{18}+
13\cdot2^{14}+17\cdot2^{11}+13\cdot2^{10}\\
&&+13\cdot2^6+17\cdot2^3+13\cdot2^2.
\end{eqnarray*}
Note that $\tS(M)=\{13,17\}$. By Corollary \ref{C:2} we have
\fbox{$c_{13}=0$}.

Suppose $c_{13}=0$. Let $\supp X=\{3,5,9,11,17\}$. We checked that
for $\lambda=\frac{1}{2}$, $j=4$, $n_0=19$ and $m_0=1$ the conditions
of Lemma~\ref{L:key} b) are satisfied, and
$\tK(m_02^{n_0+g-1}-4 )$ consists of a unique element $M$ as
follows
\begin{eqnarray*}
2^{26}-4=11\cdot2^{22}+17\cdot2^{20}+11\cdot2^{18}+
11\cdot2^{14}+17\cdot2^{12}+11\cdot2^{10}+\\11\cdot2^6+17\cdot2^4+11\cdot2^2.
\end{eqnarray*}
Note that $\tS(M)=\{11,17\}$. By Corollary \ref{C:2} we have
\fbox{$c_{11}=0$}.

\subsection{$g=9$}

Suppose $X$ is a  hyperelliptic supersingular curve over
$F$ of genus $9$. By Theorem~\ref{MT:3} II and
Proposition~\ref{P:1}
it has an equation
$$\fbox{$y^2-y=x^{19}+c_{13}x^{13}+c_{11}x^{11}
+c_9x^9+c_7x^7+c_5x^5+c_3x^3+c_1x$}$$
for some $c_\ell\in F$.

Let $\supp X=\{1,3,5,7,9,11,13,19\}$. We found for
$k=0,\ldots,7$ that
$\tK(2^{27-3k}-8)$ consists of
all sequences $\{[b_i,\ell_i]\}_{i=1}^{8-k}$ possessing the following
property:
$$
\left\{
\begin{aligned}
&[b_{8-k},\ell_{8-k}]=[24-3k,7] \mbox{ or}\\
&[b_{8-k},\ell_{8-k}]=[23-3k,13] \mbox{ and }
  [b_{7-k},\ell_{7-k}]=[21-3k,11] \mbox{ or }\\
&[b_{8-k},\ell_{8-k}]=[23-3k,11]  \mbox{ and }
  [b_{7-k},\ell_{7-k}]=[21-3k,19]  \mbox{ or }\\
&[b_{8-k},\ell_{8-k}]=[23-3k,13]  \mbox{ and }
  [b_{7-k},\ell_{7-k}]=[20-3k,19] \mbox{ and }
  [b_{6-k},\ell_{6-k}]=[18-3k,19].
\end{aligned}
\right.
$$
Hence
\begin{eqnarray*}
\tC(2^{27-3k}-8)&=&c_7^{2^{24-3k}}\tC(2^{24-3k}-8)+(c_{11}c_{13}^4+
c_{11}^4)^{2^{21-3k}}\tC(2^{21-3k}-8)\\
&&+c_{13}^{2^{23-3k}}\tC(2^{18-3k}-8).
\end{eqnarray*}
We checked that the conditions of Lemma~\ref{L:key} b)
are satisfied for
$\lambda=\frac{1}{2}$, $j=8$  and $(n_0,m_0)$ equal to
$(17,4)$, $(15,2)$ and $(13,1)$. Using the same kind of reasoning as
for genus 8, one concludes that
\fbox{$c_7=0$}, \fbox{$c_{13}=0$} and $c_{11}c_{13}^4+c_{11}^4=0$,
which clearly implies  \fbox{$c_{11}=0$}.

So we continue with $X=\{1,3,5,9,19\}$.
We checked that for $n_0=19$, $m_0=1$, $j=8$
the conditions of Lemma~\ref{L:key} b) are satisfied, and
$\tK(m_02^{n_0+g-1}-8)$ consists of
a unique element $M$ as follows
$$
2^{27}-8=
5\cdot 2^{24}+19\cdot 2^{21}+19\cdot 2^{19}
+5\cdot 2^{16}+19\cdot 2^{13}+19\cdot 2^{11}
+5\cdot 2^{8}+19\cdot 2^{5}+19\cdot 2^{3}.
$$
Note that $\tS(M)=\{5,19\}$.
By Corollary \ref{C:2} we have \fbox{$c_5=0$}.

Now set $ \supp X=\{1,3,9,19\}$.
We checked that for $n_0=25$, $m_0=3$, $j=8$
the conditions of Lemma~\ref{L:key} b) are satisfied, and
$\tK(m_02^{n_0+g-1}-8)$ consists of
a unique element $M$ as follows
\begin{eqnarray*}
3\cdot2^{33}-8&=&
19\cdot2^{30}+19\cdot2^{28}+3\cdot2^{25}+19\cdot2^{23}+
3\cdot2^{20}+19\cdot2^{18}+3\cdot2^{15}\\
&&+19\cdot2^{13}+3\cdot2^{10}+
19\cdot2^{8}+3\cdot2^5+19\cdot2^{3}.
\end{eqnarray*}
Note that $\tS(M)=\{3,19\}$.
By Corollary \ref{C:2} we have \fbox{$c_3=0$}.

Finally, set $ \supp X=\{1,9,19\}$. We found that for $n_0=35$,
$m_0=3$, $j=8$ and $\lambda=\frac{1}{2}$ the conditions of
Lemma~\ref{L:key} b) are satisfied, and that
$\tK(3\cdot2^{43}-8)$ consists of sequences
$\{[b_i,\ell_i]\}_{i=1}^{17}$ with
$[b_{17},\ell_{17}]=[40,19]$, $[b_{16},\ell_{16}]=[38,19]$
and for $k=1,\ldots,5$ we either have
$[b_{3k},\ell_{3k}]=[7k-4,19]$, $[b_{3k+1},\ell_{3k+1}]=[7k-2,19]$,
$[b_{3k+2},\ell_{3k+2}]=[7k+1,1]$ or
$[b_{3k},\ell_{3k}]=[7k-4,19]$, $[b_{3k+1},\ell_{3k+1}]=[7k-2,9]$,
$[b_{3k+2},\ell_{3k+2}]=[7k-1,9]$.
So from Lemma~\ref{L:key} it follows that
\begin{eqnarray*}
0=\tC(3\cdot2^{43-8}-8)=(c_1^{32}+c_9^{12})^
{2^3+2^{10}+2^{17}+2^{24}+2^{31}}.
\end{eqnarray*}
Thus \fbox{$c_9^3=c_1^8$}, and we conclude that
every hyperelliptic supersingular curve of genus $9$
over $F$ has to have an equation of the form
$y^2-y=x^{19}+c^8x^9+c^3x$ for some $c\in F$.

A straightforward computation shows that the curve
$X$ with equation $y^2-y=x^{19}+c^8x^9+c^3x$ over $\overline{\F}_2$ is
supersingular for $c\in\F_2$ but is {\em not} supersingular for
$c\in\F_{2^2}-\F_2$. Therefore, this curve is supersingular for only
finitely many $c\in\overline\F_2$. In other words, the open locus
of hyperelliptic supersingular curves of genus $9$ over $\overline\F_2$
is of dimension $0$.

Let $c\in F-\overline\F_2$. Supose $X:y^2-y=x^{19}+c^8x^9+c^3x$ is
supersingular. Then its specialization yields infinitely many
supersingular curves over $\overline\F_2$ because the Grothendieck
Specialization theorem says that the Newton polygon goes up under
specialization maps (see \cite{Katz}). This would contradicts
the conclusion in the previous paragraph. Therefore,
$X$ has to be defined over $\overline\F_2$.

\begin{question}
What are normal forms of these hyperelliptic supersingular curve of
genus $9$ over $\overline\F_2$?  Are $y^2-y=x^{19}$ and
$y^2-y=x^{19}+x^{9}+x$ the only ones?
\end{question}

Computations by our algorithm reveal some pattern of the dimension of
the open locus $\HS_g/\overline\F_2$ of hyperelliptic supersingular
curves of genus $g\geq 10$ over $\overline\F_2$. Vaguely speaking
$\dim(\HS_g/\overline\F_2)$ is related to $s(2g+1)$ in the sense that
$\dim(\HS_g/\overline\F_2)$ goes smaller as the ratio of $1$'s in the
binary expansion of $2g+1$ goes bigger. This is already evident in our
previous paper \cite{SZ:1} which shows that for $2g+1= 2^h-1$ for some
$h>2$ then $\HS_g/\overline\F_2$ is empty. This paper further explores
this relationship.

\begin{question}
It is desirable to formulate the pattern of $\dim(\HS_g/\overline\F_2)$
more precisely
by more numerical data.
\end{question}

\end{document}